\documentclass{amsart}
\usepackage{amssymb,psfrag,graphicx}
\newcommand{\N}{\mathbb{N}}
\newcommand{\R}{\mathbb{R}}
\newcommand{\Z}{\mathbb{Z}}
\newcommand{\norm}[3]{\ensuremath{\left\Vert#1\right\Vert_{#2}^{#3}}}
\newtheorem{bigthm}{Theorem}
\newtheorem{thm}{Theorem}[section]
\newtheorem{prop}[thm]{Proposition}
\newtheorem{cor}[thm]{Corollary}
\newtheorem{quest}[thm]{Question}
\newtheorem{lem}[thm]{Lemma}
\newtheorem{fact}[thm]{Fact}
\newtheorem*{cor*}{Corollary}

\theoremstyle{remark}
\newtheorem{rem}[thm]{Remark}
\newtheorem*{rem*}{Remark}

\theoremstyle{definition}

\newtheorem*{defn*}{Definition}

\newcommand{\intpart}[1]{\lfloor #1\rfloor}

\begin{document}
\title{Rates of divergence of nonconventional ergodic averages}
\author{Anthony Quas}
\address{Department of Mathematical Sciences, University of Memphis,
  Memphis, TN 38152, USA}
\email{aquas@memphis.edu}
\author{M\'at\'e Wierdl}

\email{mw@csi.hu}
\thanks{This research was supported by research grants from the NSF}

\begin{abstract}
  We study the rate of growth of ergodic sums along a sequence $(a_n)$
  of times: $S_Nf(x)=\sum_{n\le N}f(T^{a_n}x)$. We characterize the
  maximal rate of growth and identify a number of sequences such as
  $a_n=2^n$, along which the maximal rate of growth is achieved.
  
  We also return to Khintchine's Strong Uniform Distribution
  Conjecture that the averages $(1/N)\sum_{n\le N}f(nx\bmod 1)$
  converge pointwise to $\int f$ for integrable functions $f$, giving
  an elementary counterexample and proving that divergence occurs at
  the maximal rate.

\end{abstract}

\maketitle

\tableofcontents

\section{Introduction}

As the statements of many of the results below are fairly technical,
we start by stating some results that are formally corollaries of the
main theorems, but which in fact motivated the work in the paper.

We consider ergodic sums along a sequence $(a_n)$:
$S_Nf(x)=\sum_{n=1}^N f(T^{a_n}x)$ and ask for the maximal growth rate
of these sums.

\begin{bigthm}\ \label{bigthm:L1}
  \begin{enumerate}
  \item Let $f\in L^1$, let $T$ be a measure-preserving transformation 
    and let $(a_n)$ be an arbitrary sequence. Then
    \begin{equation*}
      \lim_{N\to\infty}\frac{1}{N\log N\log\log
      N\ldots}\sum_{n=1}^Nf(T^{a_n}x)=
      0\text{ a.e.},
    \end{equation*}
    where the product in the denominator is taken over those terms
    that exceed 1.
  \item Let $T$ be an aperiodic measure-preserving transformation and
  let the sequence $(M_N)$ satisfy $M_N/(N\log N\log\log N\ldots)\to
  0$. Then there exists an $f\in L^1$ and a sequence $(a_n)$ such that
    \begin{equation*}
      \limsup_{N\to\infty}\frac{1}{M_N}\sum_{n=1}^Nf(T^{a_n}x)=\infty
      \text{ a.e.}.
    \end{equation*}
  Further, the sequence $(a_n)$ may be taken to be the sequence $(2^n)$.
  \end{enumerate}
\end{bigthm}

\begin{rem*}
  In \cite{AJR}, Akcoglu, Jones, and Rosenblatt proved that if
  $\sum_{N=1}^\infty 1/M_N$ is finite, then $(1/M_N)S_Nf$ is
  convergent for $f\in L^1$ and also it was demonstrated that if $M_N$
  is taken to be any sequence of the form $N\log N\ldots\log^{(k)}N$
  (where $\log^{(k)}$ denotes the $k$-fold composition of $\log$),
  then there exists $f\in L^1$ for which $(1/M_N)S_Nf$ is divergent.
  Based on this, they conjectured that $(1/M_N)S_Nf$ is convergent if
  and only if $\sum_N 1/M_N$ is finite. However, one can check that
  the example $M_N=N\log N\log\log N\ldots$ disproves this conjecture.
\end{rem*}

\begin{bigthm}Let $p>1$. \label{bigthm:Lp}
  \begin{enumerate}
  \item Let $f\in L^p$, let $T$ be a measure-preserving transformation 
    and let $(a_n)$ be an arbitrary sequence. Then
    \begin{equation*}
      \lim_{N\to\infty}\frac{1}{N(\log N)^{1/p}}\sum_{n=1}^Nf(T^{a_n}x)=
      0\text{ a.e.}.
    \end{equation*}
  \item Let $T$ be an aperiodic measure-preserving transformation and
    let the sequence $(M_N)$ satisfy $M_N/(N(\log N)^{1/p})\to 0$. Then
    there exists an $f\in L^p$ and a sequence $(a_n)$ such that
    \begin{equation*}
      \limsup_{N\to\infty}\frac{1}{M_N}\sum_{n=1}^Nf(T^{a_n}x)=\infty
      \text{ a.e.}.
    \end{equation*}
  Further, the sequence $(a_n)$ may be taken to be the sequence $(2^n)$.
  \end{enumerate}
\end{bigthm}

In another aspect of the paper, we consider the averages introduced by
Khintchine: For $f\in L^1([0,1))$, $K_Nf(x)=\sum_{n\le N}f(nx\bmod
1)$. Khintchine \cite{Khintchine} in 1923 conjectured that $(1/N)K_Nf(x)$
converges to the integral of $f$. This was shown to be false by
Marstrand \cite{Marstrand} in 1970. Later, Bourgain \cite{BourgainEntropy}
gave an alternative proof using his entropy method. In Section
\ref{sect:Khintch}, we give a very simple and brief demonstration of
Marstrand's result using Rokhlin towers. In fact, we show more: we
demonstrate that for suitable $f\in L^p$, the growth rate of $K_Nf$ is
exactly the same as the maximal growth rate obtained in Theorem
\ref{bigthm:Lp}.

These techniques also allow us to resolve a question of Nair
\cite{NairI} concerning the Khinchine averages taken along a
multiplicative subsemigroup of the natural numbers rather than all
natural numbers.  Our results demonstrate that the averages $(1/|G\cap
[1,N]|)\sum_{\{n\in G\colon n\le N\}}f(nx\bmod 1)$ converges for $f\in
L^1$ to the integral if and only if the semigroup $G$ is a
subsemigroup of one that is finitely generated.

We would like to thank Ciprian Demeter for making available to us his
preprint \cite{Demeter}. Many of the ideas in that paper were crucial
to us in formulating the results of Section \ref{sect:ub}. We would
also like to thank Roger Jones and Joe Rosenblatt for stimulating
discussions.

\section{Background and Statement of Results}

We will make extensive use in what follows of so-called weak
$(L^{p,\infty})$ norms. Given a function $f$ on a measure space
$(X,\mu)$, its $L^{p,\infty}$ norm is defined by
\begin{equation*}
  \|f\|_{p,\infty}=\sup_y y\cdot\mu\{x\colon |f(x)|\ge y\}^{1/p}.
\end{equation*}
In the case of a sequence $(w_t)$, its norm analagously is
$\|w\|_{p,\infty}=\sup_y y\cdot |\{t\colon w_t\ge y\}|^{1/p}$.  As is
well-known, these ``norms'' fail to be sub-additive. In the case
$p>1$, there is a true norm $|||\cdot|||_{p,\infty}$ and a constant
$C>1$ such that $\|\cdot\|_{p,\infty}/C \le |||\cdot|||_{p,\infty} \le
C\|\cdot\|_{p,\infty}$. In the case $p=1$ however, there is no
equivalent norm. For more details about these norms, the reader is
referred to Bennett and Sharpley's book \cite{BennettSharpley}.

In this paper, we will consider almost everywhere convergence of
sequences of the form $w_tA_tf(x)$, where $(w_t)$ is a sequence of
real numbers and the $A_t$ are averaging operators of various kinds.
The typical example that we will consider is the case
$A_tf(x)=(1/2^t)\sum_{n\le 2^t}f(T^{a_n}x)$, where $(a_n)$ is a
sequence of times and $T$ is a measure-preserving transformation.  A
key tool in our work will be the maximal operator $(wA)^*f(x)=\sup_t
w_tA_tf(x)$. We will say that the sequence of operators $(w_tA_t)$
\emph{satisfies a weak $(p,p)$ maximal inequality} if there exists a
constant $C>0$ such that $\|(wA)^*f\|_{p,\infty}\le C\|f\|_p$ for all
$f\in L^p$.

\begin{fact}\label{fact:div}
  Under conditions that are satisfied by all of the operators that we
  consider in this paper, we have the following:
  \begin{enumerate}
  \item (Banach Principle) If the sequence $(w_tA_t)$ satisfies a weak
    $(p,p)$ maximal inequality, then the set of functions $f\in L^p$
    for which $w_tA_tf(x)$ is convergent almost everywhere is a closed
    set in $L^p$.
  \item If the sequence $(w_tA_t)$ fails to satisfy a weak $(p,p)$
    maximal inequality in one measure-preserving system, then there is
    a function $f\in L^p$ such that
    $\limsup_{t\to\infty}w_tA_tf(x)=\infty$ almost everywhere in any
    of the systems that we consider.
  \end{enumerate}
\end{fact}

The first statement is well known (see for example Rosenblatt and
Wierdl's article \cite{RosWie} or Garsia's book \cite{Garsia}), holding
under very mild conditions on the operator. Since in this paper,
convergence will hold trivially on the dense set of bounded measurable
functions, in order to prove a positive result, it will be
sufficient to establish a maximal inequality

The second statement is based on the folklore transference principle
(see for example \cite{RosWie}), a theorem of Sawyer \cite{Sawyer} and
an adaptation appearing in an article of Akcoglu, Bellow, Jones,
Losert, Reinhold-Larsson and Wierdl \cite{ABJLRW}.  First, the
transference principle tells us that if a maximal inequality fails in
one measure-preserving system (or flow) along some given sequence of
times, then the maximal inequality fails in all measure-preserving
systems (or flows) along the sequence of times.  Sawyer proves that if
a sequence of operators on a finite measure space fails to satisfy a
maximal inequality and commutes with a ``mixing family'' of
transformations, then there exists a function giving divergence almost
everywhere. The paper \cite{ABJLRW} reaches the same conclusion for a
family of operators that are averages of iterates of a single
aperiodic measure-preserving transformation.

We will use Iverson notation for indicator functions so that by the
expression $[y<w_t<2^ty]$, we will mean the function that is equal to
1 when the condition is satisfied and 0 otherwise.  For a sequence
$(w_t)_{t\in\mathbb{N}}$ of positive real numbers, define
\begin{equation*}
  C_p(w)=
  \begin{cases} \sup_y \sum_{t} [y<w_t<2^ty]w_t
    &\text{for $p=1$}\\
    \|w\|_{p,\infty}&\text{for $1<p<\infty$}.
  \end{cases}
\end{equation*}

We note that it is convenient to formulate the results not in terms of
the ergodic sums up to $N$ as was done in the introduction, but rather
to consider the ergodic sums (or equivalently ergodic averages) up to
$2^t$. We justify this restriction as follows. First, we observe that
it is sufficient to establish convergence to 0 for non-negative
functions.  Let $(u_n)$ be a sequence of real numbers and let
$v_t=\max_{2^{t-1}<n\le 2^t}u_n$. We will show that the following
three statements are equivalent:
\begin{enumerate}
\item $u_N\sum_{n\le N}f(T^{a_n}x)\to0$ a.e. $x$, for all
  $f\in L^p$ and every sequence $(a_n)$.\label{cond:1}
\item $v_t\sum_{n\le 2^{t-1}}f(T^{a_n}x)\to0$ a.e. $x$, for all $f\in
  L^p$ and every sequence $(a_n)$.\label{cond:2}
\item $v_t\sum_{n\le 2^{t}}f(T^{a_n}x)\to0$ a.e. $x$, for all $f\in
  L^p$ and every sequence $(a_n)$.\label{cond:3}
\end{enumerate}
To see this, note that for any $t$, $x$ and $N$ satisfying
$2^{t-1}<N\le 2^t$,
\begin{equation*}
  v_t\sum_{n\le 2^{t-1}}f(T^{a_n}x)\le 2u_N\sum_{n\le
    N}f(T^{a_n}x)\le 4v_t\sum_{n\le 2^{t}}f(T^{a_n}x).
\end{equation*}
It follows that \eqref{cond:3} implies \eqref{cond:1} implies \eqref{cond:2}.

Suppose finally that \eqref{cond:2} is satisfied. Let $(v_t)$, $(a_n)$
 and $f\in L^p$ be given. Let $b_n=a_{2n-1}$ and
 $b'_n=a_{2n}$. Applying \eqref{cond:2} separately to the sequences
 $(b_n)$ and $(b'_n)$ and summing, we deduce \eqref{cond:3}.

\begin{thm}
  \label{thm:finite1}
  Let $(w_t)$ be a sequence of positive real numbers such that
  $C_1(w)<\infty$.  For each $t\in\mathbb N$, let the set
  $\mathcal T_t$ contain at most $2^t$ measure-preserving
  transformations.  Then for $f\in L^1$,
  \begin{equation*}
    \lim_{t\to \infty} w_t2^{-t}\sum_{T\in \mathcal
      T_t} f(Tx) = 0 
  \end{equation*}
  almost everywhere.
\end{thm}

\begin{thm}
  \label{thm:finite2}
  Let $1<r<p<\infty$ and let $(w_t)$ be a sequence of positive real
  numbers such that $\|w\|_{p,\infty}<\infty$. For each $t\in\mathbb N$,
  let $A_t$ be an $L^r-L^\infty$ contraction. Then for any $f\in L^p$,
  \begin{equation*}
    \lim_{t\to \infty} w_tA_tf(x)=0
  \end{equation*}
  almost everywhere.
  
  In particular, if for each $t\in\mathbb N$, the set $\mathcal
  T_t$ contains at most $2^t$ measure-preserving transformations.  Then
  for $f\in L^p$,
  \begin{equation*}
    \lim_{t\to \infty} w_t2^{-t} \sum_{T\in \mathcal T_t} f(Tx) = 0 
  \end{equation*}
  almost everywhere.
\end{thm}

\begin{thm}\label{thm:negative}
  Let $1\le p<\infty$ and let $(w_t)$ be a sequence of positive
  real numbers such that $C_p(w)=\infty$.  Let $T$ be an aperiodic
  probability measure-preserving transformation. Then there is a
  sequence $(a_n)$ of integers so that the maximal function of the
  averages
  \begin{equation*}
    w_t2^{-t}\sum_{n\le 2^t} f(T^{a_n} x)
  \end{equation*}
  is not weak $(p,p)$ and hence there exists an $f\in L^p$ for which
  the averages diverge almost everywhere.
\end{thm}

The following proposition gives a simple description of the $w_t$ for
which $C_p(w_t)<\infty$ in the case that the $w_t$ are a sufficiently
regularly decaying sequence.

\begin{prop}\label{prop:Hardy'}
    Let $(w_t)$ be a sequence of weights and let $\Phi(t)=t\log
  t\log\log t\ldots$ be defined to be the product of $t$ and all
  iterates of $\log$ that are defined and greater than 1 at $t$.
  Let $1<p<\infty$.
  \begin{enumerate}
    \item If there exists a $K$ such that $w(t)\le K/\Phi(t)$, then
      $C_1(w)<\infty$.
    \item If $w(t)\Phi(t)\to\infty$ as $t\to\infty$, then
    $C_1(w)=\infty$.
    \item If there exists a $K$ such that $w(t)\le Kt^{-1/p}$, then
    $C_p(w)<\infty$. 
    \item If $w(t)t^{1/p}\to\infty$ as $t\to\infty$, then
    $C_p(w)=\infty$.
  \end{enumerate}
\end{prop}

\begin{rem}
  \ 
  \begin{itemize}
  \item Notice that Theorems \ref{thm:finite1}, \ref{thm:finite2} and
    \ref{thm:negative} give a dichotomy: in any $L^p$, $(p\ge1)$, if
    $C_p(w)$ is finite then the averages
    $w_t2^{-t}\sum_{n<2^t} f(T^{a_n}x)$ converge along all
    sequences of times $(a_n)$ for all $L^p$ functions $f$, whereas if
    $C_p(w)$ is infinite then in every aperiodic dynamical system
    there exists a sequence of times $(a_n)$ and an $L^p$ function $f$
    for which the averages fail to converge.
  \item We strengthen this dichotomy below by showing that there are
    sequences of times $(a_n)$ that can be chosen independently of $w$
    such that if $C_p(w)=\infty$, then in every aperiodic dynamical
    system, there exists an $f$ in $L^p$ such that
    $w_t2^{-t}\sum_{n<2^t} f(T^{a_n}x)$ diverges almost
    everywhere.
  \item Theorem \ref{thm:finite1} fails in the $L^1$ case if the
    transformations are taken to be $L^1-L^\infty$ contractions. Also
    Theorem \ref{thm:finite2} fails in the $L^p$ case if the
    contractions are only assumed to be $L^p-L^\infty$ contractions.
  \end{itemize}
\end{rem}

\section{Proofs of maximal rate theorems}

We make the following observations concerning the relationships of
$L^p$ goodness for various $L^p$:

If $1<p,q<\infty$, then $C_p(w)<\infty$ if and only if
$C_q(w^{p/q})<\infty$. To see this, note that $C_p(w)<\infty$ if and
only if $\|w^p\|_{1,\infty}<\infty$ if and only if
$\|w^{p/q}\|_{q,\infty}<\infty$.

If $C_1(w)<\infty$ then $C_p(w^{1/p})<\infty$ for $p>1$. To
see this, we note that $w_n\le C_1(w)$ for all $n$ and argue as follows:
\begin{align*}
  y\#\{n\colon w_n>y\}&=y\#\{n\colon w_n\ge 2^ny\}+y\#\{n\colon
  y<w_n<2^ny\}\\
  &\le y\#\{n\colon 2^n\le C_1(w)/y\}+\sum_n[y<w_n<2^ny]w_n\\
  &\le 2C_1(w).
\end{align*}
This shows that $\|w\|_{1,\infty}\le 2C_1(w)$ so that
$C_p(w^{1/p})=\|w^{1/p}\|_{p,\infty}=\|w\|_{1,\infty}^{1/p}<\infty$.
The converse to this assertion fails as is seen by considering
$w_t=1/t$.

\begin{proof}[Proof of Proposition \ref{prop:Hardy'}]
  We deal first with the equivalence $\Phi(t)w_t$ is bounded above if
  and only if $C_1(w)<\infty$.
  
  We start by defining a quantity $C_1'(w)$ such that $C_1(w)=\infty$
  if and only if $C_1'(w)=\infty$. Namely, define
  \begin{equation*}
    C_1'(w)=\sup_z \sum_t [t>z\text{ and }w_t>2^{-z}]w_t.
  \end{equation*}
  Writing $y=2^{-z}$, this may be rewritten $C_1'(w)=\sup_y\sum_t
  [y<w_t\text{ and }1<2^ty]w_t$. 
  Comparing with $C_1(w)=\sup_y\sum_t [y<w_t<2^ty]w_t$, we see that
  \begin{equation*}
    |C_1'(w)-C_1(w)|\le \sum_t [w_t<2^ty\le 1]w_t + \sum_t[w_t\ge
    2^ty>1]w_t.
  \end{equation*}
  If $\limsup w_t>0$, it is easy to see that both $C_1(w)$ and
  $C_1'(w)$ are infinite. Otherwise, since there are only finitely
  many terms with $w_t>1$, the second term in the above inequality is
  finite. The first term is bounded above by $\sum_t [2^ty\le
  1]2^ty\le 2$ showing that $|C_1'(w)-C_1(w)|<\infty$ as required.
  
  Since $\Phi(t)/t\to\infty$ and $\Phi(t/(\log t)^2)/t\to 0$, we see
  that for large $y$, $2^y/y^2<\Phi^{-1}(2^y)<2^y$.

  We consider $\sum_{y<t<\Phi^{-1}(2^y)}1/\Phi(t)$.
  A calculation by comparison with the integral shows that 
  \begin{equation}\label{eq:est}
  \begin{aligned}
      \sum_{y<t<\Phi^{-1}(2^y)}\frac{1}{\Phi(t)}&
      \sim \sum_{y<t<\Phi^{-1}(2^y)} \frac{1}{t\log t\log\log t\ldots}\\
      &\sim\int_{y}^{\Phi^{-1}(2^y)} \frac{1}{t\log t\log\log t\ldots}\,dt\\
      &\sim 1.
    \end{aligned}
  \end{equation}

  If $\Phi(t)w_t\to\infty$, we see that for large $t$,
  $1/\Phi(t)>2^{-y}$ implies $w_t>2^{-y}$ so
  \begin{align*}
    \sum_{\{t\colon t>y \text{ and }w_t>2^{-y}\}} w_t
    &\ge \sum_{\{t\colon t>y\text{ and }\Phi(t)<2^y\}} w_t\\
    &=\sum_{y<t<\Phi^{-1}(2^y)}\frac{1}{\Phi(t)}w_t\Phi(t).
  \end{align*}
  From equation \eqref{eq:est}, we see that this is divergent
  establishing part (2) of the proposition.
  
  If on the other hand, $w_t\Phi(t)$ is bounded above, we have $w_t\le
  k/\Phi(t)$. The above calculation then shows that $C_1'(w)<\infty$
  establishing part (1).
  
  For the $L^p$ case, Suppose that $t^{1/p}w_t\to\infty$.  Given any $M$,
  for $t$ greater than some $t_0$, $w_t>M/t^{1/p}$.  It follows that for
  sufficiently small $s$, the number of $t$ such that $w_t>s$ is at
  least $M^p/s^p$. Since $M$ is arbitrary, it follows that
  $C_p(w)=\infty$. This proves part (4).
  If on the other hand, $t^{1/p}w_t$ is bounded above, we have $w_t<k/t^{1/p}$
  for some $k$ so that the number of solutions to $w_t>s$ is bounded
  above by $k^p/s^p$ for all $y$ so that $C_p(w)<k$ giving condition (3).
\end{proof}

In order to prove Theorem \ref{thm:finite1}, we start with a lemma.

\begin{lem}
  \label{lem:fa1}
  Let $(w_t)$ be a sequence of positive numbers, let $\mathcal T_t$
  be a set of at most $2^t$ measure-preserving transformations of a
  probability space $X$ and denote by $A_tf(x)$, the quantity
  $1/2^t\sum_{T\in \mathcal T_t}f(Tx)$.
  Then for any $f\in L^1$, we have
  \begin{equation*}
    \|\sup_t w_tA_tf\|_{1,\infty}\le 9C_1(w)\|f\|_1
  \end{equation*}
\end{lem}

\begin{proof}
  We want to show that for every $\lambda >0$ and $f\in
  L^1$,
  \begin{equation*}
    \mu\left(\sup_t w_tA_tf(x)>\lambda\right)\le 9C_1(w){\lambda} \|f\|_1.
  \end{equation*}
  We will in fact prove the following apparently stronger
  inequality.
  \begin{equation*}
    \sum_t \mu\left ( w_tA_tf(x) > \lambda \right) \le 
    \frac{9C_1(w)}{\lambda} \|f\|_1.
  \end{equation*}
  Fix an $f\in L^1$. By rescaling $f$ if necessary, we can assume that
  $\lambda = 3$.  For a fixed $t$, let us split $f$ into three parts,
  up, middle, and down: $f=u + m + d$ where
  \begin{align*}
    u &= u_t = [ f \ge 2^t/w_t ] f\\
    m &= m_t = [ 1/w_t< f < 2^t/w_t ] f\\
    d &= d_t = [ f\le 1/w_t ] f
  \end{align*}
  
  We first estimate the upper part, $u=u_t$.  We note that the set of
  $x$ where $w_tA_tu(x)>1$ is a subset of the set of $x$ where
  $A_tu(x)>0$. The set on which $A_tu(x)>0$ is of measure at most
  $2^t$ times the measure of the set on which $u$ is supported. It
  follows that $\mu\{x\colon w_tA_tu(x)>1\}\le 2^t\mu\{x\colon
  f(x)\ge 2^t/w_t\}$.
  
  We can check that $C_1(w)<\infty$ implies that the $w_t$ are bounded
  above by $C_1(w)$. Hence summing over $t$, we get
  \begin{align*}
    &\sum_t \mu\left (w_tA_tu_t(x) > 1\right)
    \le\sum_t 2^t  \mu ( f \ge 2^t/w_t )\\
    &\le\sum_t 2^t \mu (C_1(w)f\ge 2^t)\le 2C_1(w)\|f\|_1.
  \end{align*}

  Now for the middle part, $m_t$,
  \begin{align*}
    \mu\left (w_tA_tm_t(x) > 1 \right) 
    &\le \int w_tA_tm_t= w_t\int m_t \\
    &= w_t \int [1/w_t< f < 2^t/w_t ] f
  \end{align*}
  Summing over $t$, and interchanging the summation and integration
  \begin{align*}
    &\sum_t \mu\left (w_tA_tm_t(x) > 1 \right) 
    \le \sum_t w_t \int  [1/w_t< f < 2^t/w_t ]  f   \\
    &= \int f(x) \sum_t [1/f(x) < w_t < 2^t/f(x)]w_t\,d\mu(x).
  \end{align*}
  Using the assumption that $C_1(w)<\infty$ (with the $y$ taken to be
  $1/f(x)$), we get
  \begin{equation*}
    \sum_t \mu\left (w_tA_tm_t(x) > 1 \right)\le
    \int C_1(w)f=C_1(w)\|f\|_1.
  \end{equation*}

  As for the down part, clearly, for every $t$,
  \begin{equation*}
    \left \{w_t\frac1{2^t} \sum_{T\in \mathcal T_t} d_t(Tx)
      > 1 \right\} 
    = \varnothing,
  \end{equation*}
  since $d_t \le 1/w_t$.

  Summing, we see that $\mu(\sup_t w_tA_tf>3)<3C_1(w)\|f\|_1$ as required.
\end{proof}

\begin{proof}[Proof of Theorem \ref{thm:finite1}]
  The above lemma establishes a maximal inequality. Since there is a
  dense class of bounded functions on which there is almost everywhere
  convergence, it follows that there is convergence for all $f\in L^1$
  as required.
\end{proof}

Before proving Theorem \ref{thm:finite2}, we prove a general lemma
that will imply the theorem almost immediately.

\begin{lem}\label{lem:fap}
  Let $1\leq r<p<\infty$, let $(w_t)$ be a sequence of positive
  real numbers and $A_t$ be a sequence of positive $L^r$-$L^\infty$
  contractions. Then there is a constant $C$ depending on $p$ and $r$
  such that
  \begin{equation*}
    \left\|\sup_t w_tA_t f\right\|_{p,\infty}\le
    C\|f\|_p\|w\|_{p,\infty}.
  \end{equation*}
\end{lem}

\begin{proof}
  We need to estimate $\sup_\lambda \lambda^p\mu\{\sup_t
  w_tA_tf>\lambda\}$. Since the inequality is homogeneous in $f$, it
  is sufficient to prove the estimate in the case that $\lambda=2$.
  For a fixed $n$, we write $f$ as the sum of $u_t$ and $d_t$, where
  $d_t=[f\le 1/w_t]f$ and $u_t=[f>1/w_t]f$. Since $A_t$ is an
  $L^\infty$ contraction, we see that $w_tA_td_t\le 1$ so that a
  necessary  condition for $w_tA_tf>2$ is $w_tA_tu_t>1$.
  We then estimate as
  follows:
  \begin{align*}
    &\mu\{\sup_t w_tA_tf>2\}\le\mu\{\sup_t w_tA_tu_t>1\}\\
    &\le \sum_t\mu\{A_tu_t>1/w_t\}=\sum_t\mu\{(A_tu_t)^r>w_t^{-r}\}\\
    &\le\sum_tw_t^r\int (A_tu_t)^r\,d\mu\le\sum_tw_t^r\int u_t^r\,d\mu\\
    &=\int f^r(x)\sum_t w_t^r[f(x)>1/w_t]\,d\mu(x)\\
    &=\int f^r(x)\sum_t w_t^r[w_t>1/f(x)]\,d\mu(x)
  \end{align*}
  Splitting the summation according into parts on which the $w_t$ lie
  between consecutive powers of $e^{-1}$, this is further bounded
  above by
  \begin{align*}
    & e^r\int f^r(x)\sum_{0\le j<\log f(x)}\sum_{\{t\colon e^{-j}<w_t\le
      e^{-j+1}\}}e^{-rj}\,d\mu(x)\\
    &\le e^r\int f^r(x)\sum_{0\le j<\log f(x)}e^{-rj}\#\{t\colon
    w_t>e^{-j}\}\,d\mu(x)\\
    &\le e^r\int f^r(x)\sum_{0\le j<\log
      f(x)}e^{-rj}\|w\|_{p,\infty}^pe^{pj}\,d\mu(x)\\
    &=e^r\|w\|_{p,\infty}^p\int f^r(x)\sum_{0\le j<\log f(x)}
    e^{(p-r)j}\,d\mu(x)\\
    &\le \frac{e^p}{e^{p-r}-1}\|w\|_{p,\infty}^p\int
    f^r(x)(f(x))^{p-r}\,d\mu(x)
    =\frac{e^p}{e^{p-r}-1}\|w\|_{p,\infty}^p\|f\|_p^p.
  \end{align*}
  
\end{proof}

\begin{proof}[Proof of Theorem \ref{thm:finite2}]
  The above lemma establishes an $L^p$ maximal inequality. Since for
  bounded functions there is convergence to 0 and these form a dense
  subset of $L^p$, the theorem follows.
\end{proof}

\section{Ultimate Badness}\label{sect:ub}

\begin{defn*}
  The sequence $(T_n)$ of linear operators on $L^p$ is called
  \emph{ultimately bad in $L^p$} if for any $(w_t)$ satisfying
  $C_p(w)=\infty$, the maximal function of the averages
  \begin{equation*}
    w_t\frac1{2^t} \sum_{n\le 2^t} T_{n} f(x)
  \end{equation*}
  is not weak $(p,p)$.

  A sequence $(a_n)$ of real numbers (integers) is called ultimately
  bad in $L^p$ if for any aperiodic measure-preserving flow
  $(T^t)_{t\in\R}$ (aperiodic measure-preserving transformation $T$)
  the operators $T^{a_n}$ are ultimately bad in $L^p$.
\end{defn*}

\begin{rem}\label{rem:div}
  By Fact \ref{fact:div}, in all sequences $(T_n)$ considered in this
  paper for which the above averages fail to satisfy a weak
  inequality, there is an $f\in L^p$ such that $\limsup_{t\to\infty}
  (w_t/2^t)\sum_{n\le 2^t}T_nf(x)$ is infinite almost everywhere.
\end{rem}

\begin{rem}\label{rem:subseq}
  If a sequence of transformations has a subsequence with bounded gaps
  that is ultimately bad in $L^p$, then the original sequence is also
  ultimately bad in $L^p$.
\end{rem}

We start the section by giving some equivalent formulations of
ultimate badness of sequences of times.

\begin{thm}\label{thm:BPrinc}
  Let $1\le p<\infty$. The following are equivalent.
  \begin{enumerate}
  \item \label{it:ub1}
    The sequence of times $(a_n)$ is ultimately bad
    for $L^p$.
  \item \label{it:ub2}
    There exists a $B$ such that for any sequence $(w_t)$ with
    $C_p(w)<\infty$, there is an $f\in L^p$ such that
    \begin{equation*}
      \left\|\sup_t w_tA_tf\right\|_{p,\infty}
      \ge BC_p(w)\|f\|_p,
    \end{equation*}
    where $A_tf(x)=1/2^t\sum_{j\le 2^t}f(T^{a_j}x)$.
  \end{enumerate}
\end{thm}

\begin{proof}
  Suppose we are given that condition \eqref{it:ub2} holds. Supposing
  further that $C_p(w)=\infty$, we can take truncations $w^{(n)}$ of
  $w$ with $C_p(w^{(n)})$ increasing to infinity. Then letting
  $f^{(n)}$ be the function guaranteed by the condition, we see that
  $\|\sup_t w_tA_tf^{(n)}\|_{p,\infty}\ge \|\sup_t
  w^{(n)}_tA_tf^{(n)}\|_{p,\infty}>BC_p(w^{(n)})\|f^{(n)}\|_p$.  Since
  the constants $BC_p(w^{(n)})$ increase to $\infty$, the ultimate
  badness follows so that condition \eqref{it:ub2} implies condition
  \eqref{it:ub1}.
  
  To show that condition \eqref{it:ub1} implies condition
  \eqref{it:ub2}, we argue by the contrapositive. Suppose that no
  constant $B$ as in condition \eqref{it:ub2} exists. Then for each
  $k\in\N$, there exists a sequence $(w^{(k)}_t)$ such that
  \begin{equation*}
    \sup_{\|f\|=1}\left\|\sup_t w^{(k)}_tA_tf\right\|_{p,\infty}
    \le 4^{-k}C_p(w^{(k)})
  \end{equation*}
  We may assume that the sequences $(w^{(k)}_t)$ are scaled so that
  $C_p(w^{(k)})=2^k$ and $\sup_{\|f\|=1}\left\|\sup_t
  w^{(k)}_tA_tf\right\|_{p,\infty}\le 2^{-k}$.
  Forming a new sequence $v_t=\sum_k w^{(k)}_t$, we observe that
  \begin{align*}
    \sup_{\|f\|=1}\left\|\sup_t v_tA_tf\right\|_{p,\infty}
    &=\sup_{\|f\|=1}\left\|\sup_t \sum_k w^{(k)}_t
      A_tf\right\|_{p,\infty}\\
    &\le\sup_{\|f\|=1}\left\|\sum_k \sup_t w^{(k)}_t A_t
      f\right\|_{p,\infty}.
  \end{align*}
  Since for $\|f\|_p=1$, $\|\sup_t w^{(k)}_t A_t f\|_{p,\infty}\le
  2^{-k}$, the norm of the sum is bounded above by a constant
  depending only on $p$. (In the case $p=1$, this follows from a
  result of Stein and Weiss \cite{SteinWeiss}). On the other hand,
  since $C_p(v)=\infty$, this establishes that the sequence $(a_n)$ is
  not ultimately bad in $L^p$ so that condition \eqref{it:ub1} implies
  condition \eqref{it:ub2}. 
\end{proof}

\begin{thm}\label{th:lotech}
  Let $1<p<\infty$ and let $(a_n)$ be a sequence of times.
  The following conditions are equivalent
  \begin{enumerate}
  \item The sequence of times $(a_n)$ is ultimately bad for $L^p$;
    \label{it:lotub}
  \item There exists a $C>0$ such that for all sequences of weights
    $(w_t)$ such that $\|w\|_{p,\infty}<\infty$, there exists an $f\in
    L^p$ such that $\|\sup_t w_tA_tf\|_{p,\infty}\ge
    C\|w\|_{p,\infty}\|f\|_p$. \label{it:lotub2}
  \item There exists a $C>0$ such that for all finite subsets
    $J\subset N$, there exists an $f\in L^p$ such that
    $\|\max_{j\in J}A_jf\|_{p,\infty}\ge C|J|^{1/p}\|f\|_p$.
    \label{it:lotweak}
  \item There exists a $C>0$ such that for all finite subsets
    $J\subset N$, there exists an $f\in L^p$ such that
    $\|\max_{j\in J}A_jf\|_{p}\ge C|J|^{1/p}\|f\|_p$.
    \label{it:lotstrong}
  \item The sequence of times $(a_n)$ is ultimately bad for $L^q$ for
    all $1<q<\infty$.
    \label{it:lotLq}
  \end{enumerate}
\end{thm}

The equivalence of \eqref{it:lotub} and \eqref{it:lotub2} was already
established in Theorem \ref{thm:BPrinc}.
The structure of the proof is that we first prove \eqref{it:lotub2} is
equivalent to \eqref{it:lotweak}. Since $\|f\|_p\ge \|f\|_{p,\infty}$,
we see that \eqref{it:lotweak} implies
\eqref{it:lotstrong}. Most of the work is taken up with proving the
implication \eqref{it:lotstrong} implies \eqref{it:lotweak}. The
implication \eqref{it:lotstrong} implies \eqref{it:lotLq} falls out of
the proof of this step.

\begin{rem}\label{rem:family}
  We remark that at no point in the proof do we use the fact that the
  $T^{a_n}$ are powers of a single measure-preserving
  transformation. The whole proof works verbatim if the $T^{a_n}$ are
  replaced by a family of measure-preserving transformations. We make
  use of the theorem in this form in Section \ref{sect:Khintch}.
\end{rem}

\begin{rem*}
  We note that a further by-product of the proof is the fact that in
  fact if $(a_n)$ is ultimately bad for $L^p$, then there is a $C$
  such that for every $J\subset \mathbb N$, there exists a
  characteristic function $f$ such that $\|\max_{j\in
    J}A_jf\|_{p,\infty}\ge C|J|^{1/p}\|f\|_p$.
\end{rem*}

\begin{rem*}
  In the $L^p$ case above, if we restrict to \emph{decreasing} sequences
  $(w_t)$, then condition \eqref{it:lotweak} can be weakened to
  \begin{equation*}
    \norm{\max_{1\le j \le N} A_jf }{p,\infty}{}
    > c\cdot N^{1/p}\cdot \norm{f}{p}{}.
  \end{equation*}
  This condition is the same as one appearing in recent work of
  Demeter \cite{Demeter}
\end{rem*}

\begin{proof}[Proof of Theorem \ref{th:lotech}: \eqref{it:lotub2} is
  equivalent to \eqref{it:lotweak}] 
  To see that condition \eqref{it:lotub2} implies condition
  \eqref{it:lotweak}, let $J$ be a finite subset of the positive
  integers and let $(w_t)$ be the indicator sequence of the set $J$.
  It is not hard to see that $\|w\|_{p,\infty}=|J|^{1/p}$ and condition
  \eqref{it:lotweak} follows.
  
  Now suppose that condition \eqref{it:lotweak} holds with a constant
  $C$. Let the sequence $(w_t)$ satisfy $\|w\|_{p,\infty}<\infty$.
  Let a positive number $\sigma<1$ be given and let $\lambda$ be such
  that
  \begin{equation*}
    \lambda \cdot |\{ t\colon  w_t > \lambda \}|^{1/p}
    > \sigma\cdot\norm{w_j }{p,\infty}{}.
  \end{equation*}
  Setting $J=\{j\colon w_j > \lambda \}$ the above can be written as
  \begin{equation*}
    \lambda \cdot |J|^{1/p} > \sigma\cdot\norm{w}{p,\infty}{}.
  \end{equation*}
  By condition \eqref{it:lotweak}, there exists an $f$  such that
  $\|\max_{j\in J}A_jf\|_{p,\infty}\ge C|J|^{1/p}\|f\|$.
  Now estimate as
  \begin{align*}
    \norm{\sup_j w_j\cdot A_jf }{p,\infty}{}
    & \ge \norm{ \max_{j\in J} w_j\cdot A_jf }{p,\infty}{}\\
    & > \norm{\max_{j\in J} \lambda\cdot A_jf }{p,\infty}{} \\
    &\ge C\cdot\lambda \cdot |J|^{1/p} \norm{f}{p}{} \\
    &> C\cdot\sigma \cdot\norm{w_j }{p,\infty}{} \norm{f}{p}{},
  \end{align*}
  where the third inequality comes from condition \eqref{it:lotweak}. This
  shows that condition \eqref{it:lotub2} follows.
\end{proof}

The proof of \eqref{it:lotstrong} implies \eqref{it:lotweak} proceeds
by three lemmas.

\begin{lem}\label{lem:refine}
  Suppose that for some finite subset $J$ of $\mathbb N$ and some
  function $f\in L^p$, $\|\max_{j\in J}A_jf\|_p=C|J|^{1/p}\|f\|_p$.
  Then there exists a subset $J'$ of $J$ such that for $j\in J'$,
  $\|A_jf\|\ge(C/2)\|f\|$ and $\|\max_{j\in J'}A_jf\|_p\ge
  (C/2)|J|^{1/p}\|f\|_p$.
\end{lem}

\begin{proof}
  Let $J_1=\{j\colon \|A_jf\|<(C/2)\|f\|_p\}$ and $J'=J\setminus
  J_1$.
  For $j\in J$, let $E_j=\{x\colon A_jf(x)=\max_{k\in J}A_kf(x)\}$.
  We have
  
  \begin{equation}\label{eq:junk}
    \begin{aligned}
      \sum_{j\in J'}\int_{E_j}(A_jf)^p&=\int \max_{j\in J}(A_jf(x))^p-
      \sum_{j\in J_1}\int_{E_j}(A_jf)^p\\
      &\ge C^p|J|\|f\|^p-\sum_{j\in J_1}\int (A_jf)^p\\
      &\ge C^p|J|\|f\|^p-|J|(C/2)^p\|f\|^p\ge (C/2)^p|J|\|f\|^p.
    \end{aligned}
  \end{equation}
  The conclusion then follows: $\|\max_{j\in J'}A_jf\|^p\ge \sum_{j\in
  J'}\int_{E_j}(A_jf)^p \ge (C/2)^p|J|\|f\|^p$

\end{proof}

By an \emph{averaging operator}, $A$, we will mean an operator of the
form $Af(x)=\frac 1N\sum_{n\le N}f(T^{a_n}x)$. Given an averaging
operator $A$, a fixed non-negative function $f$ and a real number
$R>1$, let $E^k$ denote $\{x\colon Af(x)\in (R^{k-1/2},R^{k+1/2}]\}$.
Given all of this, we define for $k\in \mathbb Z$, 
\begin{equation*}
  B^{k}g(x)=\frac1N\sum_{\{n\le N\colon f(T^{a_n}x)\in
    (R^{k-1},R^{k+1}]\}}g(T^{a_n}x)\mathbf 1_{E^{k}}.
\end{equation*}
Note that the range of the summation does not depend on $g$ so that $B^{k}$
is a linear operator. We modify this giving
$B'^kg(x)=B^kg(x)\mathbf{1}_{\{B^kg(x)>R^{k-1}\}}(x)$.
Also define $Bg(x)=\sum_k B^kg(x)$ and $B'g(x)=\sum_k B'^kg(x)$.

\begin{lem}\label{lem:equi}
  Let $0<C<1$ be given. There exists an $R$ with the following property:
  If $A$ is an averaging operator such that $\|Af\|_p\ge (C/2)\|f\|$,
  then if $B'$ is defined as above, we have $\|Af-B'f\|\le \|Af\|/2$.
\end{lem}
 
\begin{proof}
  Let $L$ be chosen so that $(2/C)L^{-1/p}=1/8$ and let the quantity
  $R$ in the statement of the lemma be chosen so that
  $\max(L/R^{(p-1)/2},R^{-1/2})=1/8$. Note that $R$ depends only on
  $C$.
  
  First define
  \begin{equation*}
    \rho(x)=\frac{1}{N}\sum_{n=1}^N\left(\frac{f(T^{a_n}x)}{Af(x)}\right)^p.
  \end{equation*}

  We note that $\int (Af)^p\rho=\int f^p \le
  \left(\frac{2}{C}\right)^p\int (Af)^p$.
  Given this, we estimate $\|Af-Bf\|$ in three parts.

  \begin{equation}
    \label{eq:decomp}
    \begin{aligned}
      &\|Af-Bf\|\le \left(\int_{\{x\colon
          \rho(x)>L\}}|Af(x)|^p\right)^{1/p}\\
      &+\left(\int_{\{x\colon
          \rho(x)\le L\}}\left(\sum_k 1_{E^k}(x)\frac1N\sum_{\{n\colon
            f(T^{a_n}x)\le R^{k-1}}f(T^{a_n}x)\right)^{p}\right)^{1/p}\\
      &+\left(\int_{\{x\colon
          \rho(x)\le L\}}\left(\sum_k 1_{E^k}(x)\frac1N\sum_{\{n\colon
            f(T^{a_n}x)>R^{k+1}}f(T^{a_n}x)\right)^{p}\right)^{1/p}.
    \end{aligned}
  \end{equation}
  First, for the second part of \eqref{eq:decomp}, we note that for
  any $x$, $x$ belongs to some $E^k$. We then calculate
  $\frac1N\sum_{\{n\colon f(T^{a_n}x)\le R^{k-1}\}}f(T^{a_n}x)<R^{k-1}\le
  Af(x)/\sqrt{R}\le Af(x)/8$. It follows that the contribution of the
  second term is dominated by $\|Af\|/8$.

  For the first term, we have 
  \begin{align*}
    \left(\frac{2}{C}\right)^p\int (Af)^p&\ge\int_{\{x\colon
    \rho(x)>L\}}\rho\cdot (Af)^p\\
    &\ge L\int_{\{x\colon \rho(x)>L\}}(Af)^p,
  \end{align*}
  so that the first term is dominated by $(2/C)L^{-1/p}\|Af\|=\|Af\|/8$.
  
  Finally, for the third term, let $x$ satisfy $\rho(x)\le L$. Since
  the $E_k$ partition the space, we let $x\in E_k$. We have
  \begin{align*}
    L&\ge \frac{1}{N}\sum_{\{n\colon f(T^{a_n}x)>R^{k+1}\}}
    \left(\frac{f(T^{a_n}x)}{Af(x)}\right)^p\\
    &=\frac{1}{N}\sum_{\{n\colon f(T^{a_n}x)>R^{k+1}\}}
    \frac{f(T^{a_n}x)}{Af(x)}\left(\frac{f(T^{a_n}x)}{Af(x)}\right)^{p-1}.
  \end{align*}

  It follows that
  \begin{equation*}
    \frac{1}{N}\sum_{\{n\colon
    f(T^{a_n}x)>R^{k+1}\}}f(T^{a_n}x)R^{(p-1)/2}\le LAf(x),
  \end{equation*}
  so that $\frac{1}{N}\sum_{\{n\colon
    f(T^{a_n}x)>R^{k+1}\}}f(T^{a_n}x)\le (L/R^{(p-1)/2})Af(x)$.  We
  see that the contribution from the last term is dominated by
  $(L/R^{(p-1)/2})\|Af\|\le \|Af\|/8$.

  To complete the proof, we note that $\|Bf-B'f\|\le
  \|Af/\sqrt{R}\|\le \|Af\|/8$ so that $\|Af-B'f\|\le \|Af\|/2$.
\end{proof}

\begin{lem}\label{lem:tri}
  Let $J'\subset J$ and suppose $(F_j)_{j\in J'}$ and $(G_j)_{j\in
    J'}$ satisfy $G_j\le F_j$ and $\|G_j-F_j\|\le (C/4)\|f\|_p$.
  Suppose further that $\|\max_{j\in J'}F_j\|\ge (C/2)|J|^{1/p}\|f\|$.
  Then
  \begin{equation*}
    \|\max_{j\in J'}G_j\|\ge (C/4)|J|^{1/p}\|f\|
  \end{equation*}
\end{lem}

\begin{proof}
  Let $H=\max_{j\in J'}F_j(x)-\max_{j\in J'}G_j(x)$ and let $E_j$
  be the set $\{x\colon F_j(x)=\max_{k\in J'}F_k(x)\}$. Then
  \begin{align*}
    \|H\|^p&=\int H^p=\sum_{j\in J'}\int_{E_j} H^p\\
    &=\sum_{j\in J'}\int_{E_j}(F_j-\max_{k\in J'} G_k)^p\\
    &\le \sum_{j\in J'}\int_{E_j}(F_j-G_j)^p\\
    &\le \sum_{j\in J'}(C/4)^p\|f\|^p\le(C/4)^p|J|\cdot\|f\|^p.
  \end{align*}

We then have that $\|\max_{j\in J'}G_j\|\ge \|F\|-\|H\|\ge
(C/4)|J|^{1/p}\|f\|$.
\end{proof}

We now assemble the above lemmas to complete the proof of Theorem
\ref{th:lotech}.

\begin{proof}[Proof of Theorem \ref{th:lotech}: \eqref{it:lotstrong}
  implies \eqref{it:lotweak} and \eqref{it:lotLq}]
  Recall that we are assuming that $\|\max_{j\in J}A_jf\|=C|J|^{1/p}\|f\|$.
  From Lemma \ref{lem:refine}, we can pick a subset $J'$ of $J$ such
  that for each $j\in J'$, $\|A_jf\|\ge (C/2)\|f\|$ and also such that
  $\|\max_{j\in J'}A_jf\|\ge (C/2)|J|^{1/p}\|f\|$.

  By Lemma \ref{lem:equi}, there exists an $R$ and operators $B'_j$ for
  each $j\in J'$ such that $\|A_jf-B'_jf\|\le \|A_jf\|/2$.
  Specifically, these operators were defined as follows.
  Let $E_j^k=\{x\colon A_jf(x)\in (R^{k-1/2},R^{k+1/2}]\}$ and
  define operators by
  \begin{equation*}
    B_j^kf(x)=\frac{1}{|I_j|}\sum_{\{n\in I_j\colon f(T^{a_n}x)\in
    (R^{k-1},R^{k+1}]\}} f(T^{a_n}x)\mathbf 1_{E_j^k}(x),
  \end{equation*}
  where $I_j$ is the range of indices of the $a_n$ involved in the
  $j$th average. We then define
  ${B'}_j^kf(x)=B_j^kf(x)\mathbf{1}_{\{x\colon B_j^kf(x)>R^{k-1}\}}$ and 
  ${B'}_jf(x)=\sum {B'}_j^kf$. Note that the non-zero values taken by
  ${B'}_j^k$ are in the range $(R^{k-1},R^{k+1}]$.

  Applying Lemma \ref{lem:tri} with $F_j=A_jf$, $G_j=B_j'f$,
  we deduce $\|\max_{j\in J'}B_j'f\|\ge (C/4)\|f\|$.

  Now write $f$ as a decomposition $f=\ldots+f_{-1}+f_0+f_1+f_2+\ldots$,
  where $f_k(x)=f(x)\cdot \mathbf 1_{\{x\colon R^{k-1}<f(x)\le
  R^k\}}$. Suppose that $\max_{j\in J'}B_j'f(x)\in
  (R^{k-1},R^{k}]$. We will assume that the maximum is attained
  for $\ell\in J'$. Since the maximum is in the range
  $(R^{k-1},R^k]$, it follows that $B'_\ell f(x)={B'}_\ell^kf(x)$ or
  ${B'}_\ell^{k-1}f(x)$. In particular, we have $B'_\ell
  f(x)=B'_\ell(f_{k-2}+f_{k-1}+f_k)(x)$. Setting
  $h_k=f_{k-2}+f_{k-1}+f_k$, if $x$ satisfies $\max_{j\in J'}B_j'f(x)\in
  (R^{k-1},R^{k}]$, we have shown that $\max_{j\in J'}
  B'_jf(x)=\max_{j\in J'}B_j'h_k(x)$.

  Write $F(x)=\max_{j\in J'}B_j'f(x)$ and decompose $F$ into parts
  $F_k$ where $R^{k-1}<F(x)\le R^k$. Then the above shows that $\max_jB_j'
  h_k\ge F_k$. 
  We now have 
  \begin{align*}
    (C/12)^p|J|\sum_k\|h_k\|^p&\le (C/4)^p|J|\sum_k\|f_k\|^p\\
    &=(C/4)^p|J|\|f\|^p \\
    &\le \|F\|^p=\sum_k \|F_k\|^p\\
    &\le \sum_k \|\max_j B_j'h_k\|^p
  \end{align*}

  In particular, there must exist a $k$ such that $\|h_k\|\neq 0$ and 
  $\|\max_{j\in J'} B'_jh_k\|\ge (C/12)|J|^{1/p}\|h_k\|$.

  Since $h_k$ only takes non-zero values between $R^{k-3}$ and $R^k$,
  we see from the definition of $B'$ that $\max_{j\in J'}B'_jh_k$ only
  takes non-zero values between $R^{k-4}$ and $R^k$. It follows that
  $\|\max_{j\in J}A_jh_k\|_{p,\infty}\ge \|\max_{j\in
    J'}B'_jh_k\|_{p,\infty}\ge (C/12R^4)|J|^{1/p}\|h_k\|$.
  
  Further, since both $h_k$ and $\max_{j\in J'}B'_jh_k$ take on values
  in ranges with a bounded ratio between the endpoints, it follows
  that for any $1<q<\infty$, there exists a $C'$ such that
  $\|\max_{j\in J}A_jh_k\|_{q,\infty}\ge \|\max_{j\in
    J'}B'_jh_k\|_{q,\infty}\ge C'|J|^{1/q}\|h_k\|_q$.
  Applying the equivalence \eqref{it:lotweak} implies \eqref{it:lotub}
  for $L^q$ completes the proof of the theorem.

\end{proof}

The proof of Theorem \ref{thm:negative} will depend on the following
lemmas.

\begin{lem}
  \label{lem:ubL1}
  Suppose the sequence $(a_n)$ of real numbers
  satisfies the following condition:
  \begin{quote}
    For each positive integer $M$, there exists an $n_0$ such that if
    $N\ge n_0$ and $K$ satisfies $K\le MN$ then for any sequence
    $r_1,r_2,\dots,r_{n}$ of integers, there is a positive real
    $\alpha$ so that for 
    \begin{equation*}
      \intpart{\alpha a_{N+k}}\equiv r_k \mod{K},\;  \text{ for all
      $1\le k\le n$.} 
    \end{equation*}
  \end{quote}
  Then the sequence of times $(a_n)$ is ultimately bad for $L^1$.
\end{lem}

\begin{lem}
  \label{lem:ubLp}
  Suppose the sequence $(a_n)$ of real numbers
  satisfies the following condition:
  \begin{quote}
    There exists an $n_0$ such that if
    $N\ge n_0$ and $K$ satisfies $K\le N$ then for any sequence
    $r_1,r_2,\dots,r_{n}$ of integers, there is a positive real
    $\alpha$ so that
    \begin{equation*}
      \intpart{\alpha a_{2^{N-1}+k}}\equiv r_k \mod{K},\;  \text{ for all
      $1\le k\le n$.} 
    \end{equation*}
  \end{quote}
  Then the sequence of times $(a_n)$ is ultimately bad in every $L^p$ ($p>1$).
\end{lem}

\begin{cor}\label{cor:L1ub}
  Suppose that for some fixed $\epsilon$ the sequence $(a_n)$ satisfies 
  \begin{equation*}
    \frac{a_{n+1}/a_n}{n^\epsilon} \to \infty.
  \end{equation*}

  Then $(a_n)$ is ultimately bad in  $L^1$.
\end{cor}

\begin{proof}
  This follows from Lemma \ref{lem:ubL1} using a standard lacunarity
  argument. Let the sequence $(a_n)$ be as in the statement of the
  lemma.  Using Remark \ref{rem:subseq}, we first refine to a subsequence
  $(b_n)$ of the $(a_n)$'s occurring with bounded gaps in the original
  sequence such that $(b_{n+1}/b_n)/n^2\to\infty$. It will be
  sufficient to show that $(b_n)$ is ultimately bad for $L^1$. Let $M$
  be given.  Then there exists an $N>2M$ such that $n\ge N$ implies
  $b_{n+1}/b_n>n^2>2MN$. Let $K\le MN$, we see that for any $n\geq N$,
  $b_{n+1}/b_n>2K$. 

  To finish the argument, we claim that for any sequence
  $r_0,r_1,\ldots,r_{t-1}$ of integers with $0\le r_i<K$, there exists
  an interval $I$ of length $1/(Kb_{t-1})$ such that for $\alpha\in
  I$, $\alpha b_i\bmod 1\in [r_i/K,(r_i+1)/K)$ for $0\le i<t$.

  We prove this by induction. Clearly it is true for $t=1$. Suppose
  that it holds for $t\le s$ and let $I$ be the interval of length
  $1/(Kb_{s-1})$ such that for $\alpha\in I$, $\alpha b_i\bmod 1\in
  [r_i/K,(r_i+1)/K)$ for $0\le i<s$. We see that $S=\{\beta\colon
  b_s\beta\bmod 1\in [r_s/K,(r_s+1)/K)$ is a union of intervals of
  length $1/(Kb_s)$ spaced $1/b_s$ apart. Since
  $1/b_s(1/K+1)<1/(Bb_{s-1})$ we see that $I$ contains a complete
  interval from $S$. Letting $J$ be the subinterval, the induction is
  complete.
\end{proof}

\begin{cor}\label{cor:Lpub}
  Suppose that for some fixed $\epsilon$ the sequence $(a_n)$ satisfies 
  \begin{equation*}
    \frac{a_{n+1}/a_n}{(\log n)^\epsilon} \to \infty.
  \end{equation*}

  Then $(a_n)$ is ultimately bad in $L^p$ for every $1<p<\infty$.
\end{cor}

\begin{proof}
  The proof is similar to that of Corollary \ref{cor:L1ub}. Let
  $(a_n)$ be as in the statement and suppose that $(a_{n+1}/a_n)/(\log
  n)^\epsilon \to\infty$.  First, refine the sequence to a subsequence
  $(b_n)$ with bounded gaps in the original sequence so that the new
  sequence $(b_n)$ satisfies $(b_{n+1}/b_n)/(\log n)\to\infty$.  There
  exists an $n_0$ such that $n\ge 2^{n_0-1}$ implies
  $b_{n+1}/b_n>3\log (2n)$. Now if $K\le N$, we see that for any $n>
  2^{N-1}$, $b_{n+1}/b_n>3\log (2n)>2N\ge 2K$.  The remainder of the
  argument follows exactly as in Corollary \ref{cor:L1ub}
\end{proof}

\begin{rem*} Unfortunately both in the $L^1$ and $L^p$ cases, an
  arbitrary lacunary sequence $(a_n)$ is out of reach for now.
\end{rem*}

\begin{cor*}
  The sequence $(n!)$ is ultimately bad in every $L^p$, $1\le
  p<\infty$.
\end{cor*}
\begin{proof}
  This is an immediate consequence of Corollaries \ref{cor:L1ub} and
  \ref{cor:Lpub}.
\end{proof}

\begin{proof}[Proof of Theorem \ref{thm:negative}]
  This follows from the above Corollary.
\end{proof}

\begin{cor*}
  Let the sequence $(a_n)$ be independent over the rationals.
  Then $(a_n)$ is ultimately bad in every $L^p$, $1\le p<\infty$.
\end{cor*}

\begin{proof}
  This follows since the vectors $\alpha (a_{N},
  a_{N+1},\dots,a_{N+t-1})\bmod 1$ are dense in the $t$-dimensional
  torus as $\alpha$ runs over the positive reals. Accordingly, for any
  $N$, $K$ and $r_0,\ldots,r_{t-1}$, there exists a positive real number
  $\alpha$ with the property that $\alpha (a_{N}, a_{N+1},\dots,a_{N+t-1})\bmod
  1\in \prod_{0\le i<t}[r_i/K,(r_i+1)/K)$.
\end{proof}

\begin{cor*}
  The sequence $(\sqrt n)$ is ultimately bad in every $L^p$, $1\le
  p<\infty$.
\end{cor*}

\begin{proof} The set $\{\sqrt s\colon s \text{ squarefree}\}$ is
  independent over the rationals and arranged in increasing order it
  forms a positive density subsequence of $(\sqrt n)$ so that there
  exists a fixed $k$ such that the first $2^n$ squarefree numbers are
  a subset of the first $2^{n+k}$ square roots. We then use the fact
  that if $(w_t)$ satisfies $C_p(w)=\infty$, then
  $C_p(w_{t+k})=\infty$, so letting $A_t$ be the average over the
  first $2^t$ square roots and $B_t$ be the average over the first
  $2^t$ squarefree square roots, we have the estimate
  $w_{t+k}A_{t+k}f\ge 2^k w_{t+k}B_tf$. Since $C(w_{t+k})=\infty$, the
  right hand side can be made to diverge and hence so does the left
  hand side.
\end{proof}

\begin{proof}[Proof of Lemma~\ref{lem:ubL1}]
  
  We just deny the maximal inequality on $[0,1)$. Let the positive
  integer $M$ be given and let $n_0$ be as in the statement of the
  lemma. We consider the set of $N$ such that
  \begin{equation}
    \label{eq:ycond}
    \sum_{t} [1/N<w_t<2^t/(4N)]w_t > 3M.
  \end{equation}
  We note that these terms are necessarily unbounded above as
  $\sum_t[z<w_t<2^tz]w_t$ may be bounded above by the sum of two of
  these terms.  We assert that the set of $N$ satisfying
  \eqref{eq:ycond} is unbounded above.
  This is because either $\limsup w_t>0$, in which case for all large
  enough $N$ we have $\sum_{t} [1/N<w_t<2^t/N] w_t=\infty$, or $w_t\to
  0$, in which case for $N<K$, $\sum_{t} [1/N<w_t<2^t/N]w_t$ is
  uniformly bounded above by $\sum_t [w_t>1/K]w_t$ which is the sum of
  a finite number of terms and hence is finite. Since the sums
  $\sum_{t} [1/N<w_t<2^t/(4N)]w_t$ are as noted above unbounded in $N$,
  there must exist arbitrarily large integers $N$ for which the sum exceeds
  $3M$.  Hence we may choose an $N$ satisfying \eqref{eq:ycond} such
  that $N>n_0$.
  
  If we consider the $t$'s such that $2^t\le 2N$, then we see
  $\sum_{\{t\colon 2^t\le 2N\}}[1/N<w_t<2^t/(4N)]w_t<\sum_{\{t\colon
  2^t\le 2N\}}2^t/(4N)<1$ so that $\sum_{\{t\colon
  2^t>2N\}}[1/N<w_t<2^t/(4N)]w_t >2M$. 
  It then follows that there is a finite set $U$ of the $t$'s
  satisfying $2^t>2N$ and $1/N<w_t<2^t/(4N)$ so that we still have
  $\sum_{t\in U}w_t > 2M$.
  
  Now set $K=NM$. For each $t\in U$, select
  $\lfloor Nw_t\rfloor$ different residue classes modulo $K$.
  Denote these residue classes by $R_t$. Since
  \begin{equation*}
    \sum_{t\in U}\lfloor Nw_t\rfloor > NM,
  \end{equation*}
  we can choose the $R_t$ so that their union over $t\in U$ covers all
  residue classes modulo $K$. Since $2^{t-1}>N$ for all $t\in U$, we
  can now apply the condition in the statement of the lemma to
  conclude that there is a positive $\alpha$ so that for each $t\in U$
  and $r\in R_t$ there are at least $2^{t-1}/(Nw_t)$ $n$'s between
  $2^{t-1}$ and $2^t$ with $\intpart{\alpha a_n}\equiv r\mod{K}$.
  Define the function $f$ by
  \begin{equation*}
    f(x)=
    \begin{cases}
      2N,
      &\text{if $0\le x<\frac2K$}\\
      0, & \text{otherwise}
    \end{cases}
  \end{equation*}
  Clearly, $\|f\|_1 = \frac{4}M$.
  
  We define a measure-preserving flow on $[0,1]$ by
  $T^\zeta(x)=x+\alpha \zeta/K$.
  Now let $x\in [0,1)$ be arbitrary. Then there is a $t$ and $r\in
  R_t$ with $x\in [\frac{K-r}K,\frac{K-r+1}{K})$.
  Since we have a number of $n$ between $2^{t-1}$ and $2^t$ such
  that $\intpart{\alpha a_n}\equiv r \mod{K}$, for
  these $n$, we see that $T^{a_n}(x)\in [0,2/K)$
  so that $f(T^{a_n}(x))=2N$ and hence we can
  estimate:
  \begin{align*}
    &w_t\frac1{2^t} \sum_{n\le 2^t} f(x+a_n) \\
    &\ge w_t\frac1{2^t} \sum_{\{n\colon 2^{t-1}< n \le 2^t \}}
    \bigl[\intpart{\alpha a_n}\equiv r\mod{K} \bigr]
    f(T^{a_n}(x)) \\
    &=w_t\frac1{2^t} \sum_{\{n\colon 2^{t-1}<  n \le 2^t \}}
    \bigl[ \intpart{\alpha a_n}\equiv r\mod{K} \bigr] 2N
    \\
    &\ge w_t\frac1{2^t}  2N  \frac{2^{t-1}}{Nw_t}\\
    &=1.
  \end{align*}

  We have shown that
  \begin{equation*}
    m\left\{x\colon \sup_tw_t\frac1{2^t} \sum_{n\le 2^t}
    f(T^{a_n}(x))\ge1\right\} =1 > \frac{M}{4}\|f\|_1.
  \end{equation*}
  Since $M$ is arbitrary, the required violation of the maximal
  inequality is shown.

\end{proof}

\begin{proof}[Proof of Lemma \ref{lem:ubLp}]
  We aim to establish condition \eqref{it:lotweak} of Theorem
  \ref{th:lotech} for the sequence $(a_n)$.  Let $n_0$ be as in the
  statement of the lemma and let $J$ be a finite subset of $\mathbb
  N$. If $|J|\le 2n_0$, taking $f$ to be a constant function, we have
  $\|\max_{j\in J}A_jf\|_{p,\infty}\ge (2n_0)^{-1/p}|J|^{1/p}\|f\|_p$.
  
  If on the other hand, $|J|\ge 2n_0$, then let $N=K=\lfloor
  |J|/2\rfloor$ and let $J'=\{j_1,\ldots,j_{K}\}$ be a subset of $J$
  of size $K$ consisting of elements of $J$ at least as big as
  $|J|/2$.  By assumption, there exists an $\alpha$ such that for
  $n\in \{2^{j_\ell-1}+1,2^{j_\ell-1}+2,\ldots,2^{j_\ell}\}$, $\lfloor
  \alpha a_n\rfloor\equiv \ell-1\pmod K$. Then letting $(T^t)$ be the
  flow on $[0,K)$ given by $T^t(x)=x-\alpha t$ and $f$ be the function
  $2\cdot \mathbf 1_{[0,2)}$, we see that
  $\|f\|_p=2^{1+1/p}|K|^{-1/p}$. We also have for $x\in [n,n+1)$,
  $A_{j_n}f(x)\ge 1$. It follows that $\|\max_{j\in
    J}A_jf(x)\|_{p,\infty}=1\ge 2^{-1-2/p}|J|^{1/p}\|f\|$ as required.

\end{proof}

\section{Khintchine's Conjecture}\label{sect:Khintch}

In this section, we consider the averages arising in a conjecture due to
Khintchine.  For $f\in L^p[0,1)$, Write $T_n f(x)=f(nx\bmod 1)$.
Khintchine \cite{Khintchine} conjectured in 1923 that for every $f\in
L^1$, it is the case that $1/N\sum_{n=1}^N T_nf(x)$ converges
pointwise almost everywhere to $\int f$. This was answered negatively
by Marstrand \cite{Marstrand} in 1970. This negative result was
strengthened further in Bourgain's work using his Entropy Method
\cite{BourgainEntropy}.

We start with a lemma showing the equivalence of maximal theorems for
averages of the type $\frac{1}{N}\sum_{n\le N}f(a_nx)$
for functions $f\in L^p([0,1))$ and averages of the type
$\frac{1}{N}\sum_{n\le N}g(x-\log a_n)$
for functions $g\in L^p(\R)$.

This lemma will allow us to give a very simple demonstration of
Marstrand's result and in fact to show more: that the sequence of
operators $(T_n)$ is ultimately bad in $L^p$ for $p>1$. 

We also take up a question posed by Nair in \cite{NairI} concerning a
version of Khintchine's conjecture, where the $T_nf(x)$ are averaged
along a subsequence rather than all of the integers. Nair proved that
if the sequence $(a_n)$ is the increasing enumeration of a finitely
generated multiplicative subsemigroup of the positive integers, then
for all $f\in L^1$, the averages $1/N\sum_{n=1}^N T_{a_n}f(x)$
converge for almost every $x$ to $\int f$. He asked about the case of
averaging along an infinitely generated subsemigroup of the positive
integers.

Later, Lacroix \cite{LacroixII,LacroixIII} took up this question and
claimed that there do exist infinitely generated subsemigroups of the
integers along which the above averages converge. Unfortunately, while
the arguments in his papers appear to be correct, the result seems to
be false as they rely on an incorrect statement in Krengel's book
\cite{Krengel}.

Here, using the lemma again, we clear up the situation with an
explicit dichotomy in Theorem \ref{th:dichotomy}. If $S$ is a
multiplicative subsemigroup of the positive integers, then the
averages above converge for all $f\in L^1$ if and only if $S$ is
contained in a finitely generated subgroup of the positive integers.

\begin{lem}
  \label{lem:equiv}
  Let $(a_n)$ be any sequence of positive integers. Let
  $I_1,\ldots,I_k$ be any non-empty finite subsets of $\N$.
  Denote by $A_jf(x)$ the average $1/|I_j|\sum_{n\in I_j}f(a_nx)$
  and by $B_jg(y)$ the average $1/|I_j|\sum_{n\in I_j}g(y-\log
  a_n)$.
  Then the following are equivalent:
  \begin{enumerate}
  \item There exists an $f\in L^p[0,1)$ such that $\|\max_{j\le
      k}A_jf\|_{p,\infty}>C\|f\|_p$;\label{it:mult}
  \item There exists a $g\in L^p(\R)$ such that $\|\max_{j\le
      k}B_jg\|_{p,\infty}>C\|g\|_p$.\label{it:add}
  \end{enumerate}
\end{lem}
 
\begin{figure}
  \psfrag{Y,g}{$g$ on the Rokhlin tower in $Y$} 
  \psfrag{[0,1],f}{$f$ on the Rokhlin tower in $X$}
  \scalebox{0.7}{\includegraphics{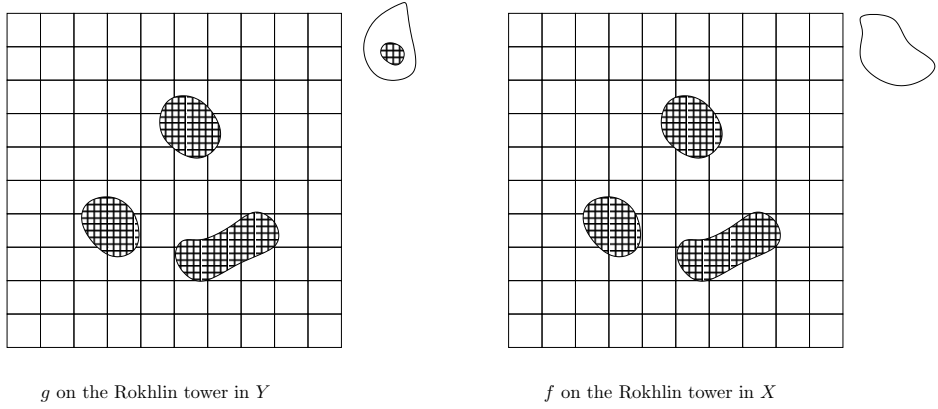}}
  \caption{Copying a function between Rokhlin towers}\label{fig:Rokhlin}
\end{figure}

Let $T_n(x)=nx\bmod 1$ and let $S_n(y)=y-\log n$. The crux of the
proof is the simple observation that $S_n$ and $T_n$ satisfy the same
basic relationship: $S_{nm}=S_n\circ S_m$ and $T_{nm}=T_n\circ T_m$,
allowing data from a Rokhlin tower for one system to be copied to a
Rokhlin tower for the other. This transference is illustrated in
Figure \ref{fig:Rokhlin}.
  
\begin{proof}
  Let $p_1,p_2,\ldots,p_d$ be the primes occurring in the prime
  factorization of elements of $\{a_n\colon n\in \bigcup I_j\}$ and let
  $r$ be the maximum of all the powers of the $p_d$ occurring in the
  elements of $\{a_n\colon n\in \bigcup I_j\}$.

  We first observe that condition \eqref{it:add} is equivalent to the
  following condition that we call (2'):
  \begin{quote}
    There exists an $M$ rationally independent of $\{\log p_1,\ldots,\log
    p_d\}$ and a $g\in L^p([0,M))$ such that 
    $\|\max_{j\le k}B_jg\|_{p,\infty}>C\|g\|_p$, where
    the difference $y-\log a_n$ is interpreted modulo $M$.
  \end{quote}
  To see that (2) implies (2'), simply restrict the function $g$
  occurring in (2) to some large interval, whereas to see that (2')
  implies (2), starting from the function in (2'), concatenate a large
  number of translated copies of the function $g$ on intervals
  $[(n-1)M,nM)$ to produce a function supported on $[0,LM)$ and
  observe that condition (2) is satisfied.
  
  We will therefore demonstrate the equivalence of (1) and (2'). If
  (2') holds, let $M$ be as in the statement, otherwise let $M=1$ so
  that $M$ is rationally independent of $\{\log p_1,\ldots,\log
  p_d\}$.  Let $N$ be chosen to be a large integer and let
  $\epsilon>0$ be small.
  
  For $\mathbf n\in\N^d$, write $T^{\mathbf n}(x)=\prod
  p_i^{n_i}x\bmod 1$, and for $\mathbf n\in \Z^d$, write $S^{\mathbf
  n}(y)=y-\sum n_i\log p_i\bmod M$. We observe that these are both free
  actions. Accordingly, we can construct Rokhlin towers of geometry
  $\Lambda_N=\{0,1,\ldots,N-1\}^d$ for both systems with an error set
  of size exactly $\epsilon$: there exist $V\subset [0,M)$ and
  $W\subset [0,1)$ such that $\mu(V)=\lambda(W)$ and that the sets
  $S^{-\mathbf n}V$ for $\mathbf n\in \Lambda_N$ are mutually
  disjoint, as are the sets  $T^{-\mathbf n}W$. We now construct
  measure-preserving maps between the Rokhlin towers. Let
  $R_X=\bigcup_{\mathbf n\in \Lambda_N}T^{-\mathbf n}W$ and
  $R_Y=\bigcup_{\mathbf n\in \Lambda_N}S^{-\mathbf n}V$. 
  
  Let $\theta_0$ be an arbitrary measure-preserving measurable
  bijection from $V$ to $W$. Then define $\theta(x)$ for $x \in
  T^{-\mathbf n}W$ by $S^{-\mathbf n}\circ \theta\circ T^{\mathbf
    n}(x)$. Similarly, letting $\mathbf k=(N-1,\ldots,N-1)$, let
  $\psi_0$ be an arbitrary measure-preserving measurable bijection
  from $S^{-\mathbf k}V$ to $T^{-\mathbf k}W$ and define $\psi(y)$ for
  $y\in S^{-\mathbf n}V$ by $\psi(y)=T^{\mathbf k-\mathbf n}\circ
  \psi_0\circ S^{-(\mathbf k-\mathbf n)}(y)$. 
  These are then defined so as to ensure that $\theta(T^\mathbf
  n(x))=S^\mathbf n(\theta (x))$ provided that the orbit of $x$
  remains inside the tower and similarly $\psi(S^\mathbf
  n(y))=T^\mathbf n(\psi(y))$.
  
  If condition (2') holds, we define $f$ on $R_X$ by
  $f(x)=g(\theta(x))$ and define $f$ to be 0 on the remainder of
  $[0,1)$. By construction, we see that provided that $x\in
  \bigcup_{\{\mathbf n\colon r\leq n_i<N, i=1,\ldots,d\}} T^{-\mathbf
    n}(W)$, we have that $f(T^\mathbf nx)=g(S^\mathbf n(\theta(x))$
  for $\mathbf n$ with coefficients less than $r$. In particular,
  since the times $a_n$ involved in the averages $A_j$ and $B_j$ for
  $j\in J$ may be expressed in terms of $p_1,\ldots,p_d$ with powers
  at most $r$, we see that for such an $x$, we have $\max_{j\in J} A_j
  f(x)=\max_{j\in J} B_j g(\theta(x))$.  Now for sufficiently small
  $\epsilon$ and large $N$, we will have $\|\max_{j\le k} A_j
  f(x)\|_{p,\infty} > C\|f\|_p$.

  If condition (1) holds, we define $g$ on $R_Y$ by $g(y)=f(\psi(y))$
  and define $g$ to be 0 elsewhere. The same argument as above
  demonstrates that condition (2') holds provided that $N$ is chosen to
  be sufficiently large and $\epsilon$ is taken to be sufficiently small.
    
\end{proof}

\begin{thm}
  \label{th:Khintchine}
  The sequence $(T_n)$ of operators defined by $T_nf(x)=f(nx\bmod 1)$
  is ultimately bad in $L^p$ for all $p>1$.
\end{thm}

\begin{proof}
  We let $g(y)$ be the function $2\cdot\mathbf 1_{[0,2\log 2)}$ and
  set for any finite set $J\subset \N$, $I_j=\{n\colon n\le 2^j\}$.
  We will then demonstrate that $\max_{j\in J}\|B_j g\|_{p,\infty}\ge
  C|J|^{1/p}\|g\|_p$ for a constant $C$ that does not depend on $J$.
  By Lemma \ref{lem:equiv}, this will establish the existence of an
  $f\in L^1[0,1)$ satisfying  condition \eqref{it:lotweak} of
  Theorem \ref{th:lotech} (see Remark \ref{rem:family}).

  We have $\|g\|_p^p=2^{p+1}\log 2$.  Let $j\in J$ and $x\in [j\log
  2,(j+1)\log 2)$. Then
  $B_j(x)\ge 1$.  It follows that the measure of the set where the
  maximal function exceeds 1 is at least $|J|$. This shows that
  $\|\max_{j\in J}B_tg(y)\|_{p,\infty}> C|J|^{1/p}\|g\|$, where
  $C=2^{-1-1/p}(\log 2)^{1/p}$ as required.
\end{proof}

\begin{lem}\label{lem:weakbound}
  Let $(h_n)$ be an non-decreasing sequence of real numbers and let
  $c_n=(h_n-h_{n-1})/h_n$ (or 0 in the case that the denominator is
  0). Then $h_n=O(n^{\|c\|_{1,\infty}})$ as $n\to\infty$.
\end{lem}

\begin{proof}
  We will suppose for simplicity that $h_1>0$. Suppose that
  $\|c\|_{1,\infty}=d<\infty$. We have $h_n=h_{n-1}/(1-c_n)$ so that
  in particular,
  \begin{equation*}
    h_n=h_1\prod_{j=2}^n \frac{1}{1-c_j}.
  \end{equation*}
  If $(t_n)$ denotes the decreasing rearrangement of $(c_n)$, then we
  have
  \begin{equation*}
    h_n\le h_1\prod_{j=2}^n \frac{1}{1-t_j}.
  \end{equation*}
  Since $\|t\|_{1,\infty}=d$, we have $|\{j\colon
  t_j>d/k\}|<d/(d/k)=k$ so that $t_k\le d/k$.
  Letting $C=h_1\prod_{j=2}^{2d-1} 1/(1-t_j)$, we have
  \begin{equation*}
    h_n\le C\prod_{j=2d}^n \frac{1}{1-d/j}.
  \end{equation*}
  Taking logarithms, we see that $\log h(n)\le C'+d\log n$ so that
  $h(n)\le Kn^d$ as required.
\end{proof}

\begin{thm}\label{thm:growth}
  Let $(t_n)_{n\in \N}$ be an increasing sequence of real numbers with
  the property that $t_n\to\infty$. Let $h(N)$ denote $|\{n\colon
  t_n\le N\}|$. If $\limsup h(N)/N^k=\infty$ for all $k$, then 
  $B_Ng$ fails to satisfy a maximal inequality, where
  $B_Ng(y)=(1/h(N))\sum_{\{n\colon t_n\le N\}}g(T_{t_n}y)$.
\end{thm}

\begin{proof}
  Set $c_n=(h(n)-h(n-1))/h(n)$. For simplicity, we assume that
  $h(1)\ge 1$. We let $g$ be the indicator function $\mathbf
  1_{[0,2)}$ and we estimate $\|\sup_N B_Ng\|$.  We quickly see
  that for $n\le x<(n+1)$, $B^*g(x)\ge
  B_{n}g(x)\ge(h(n)-h(n-1))/h(n)=c_n$. It follows that
  $\|\sup_N B_Ng\|_{1,\infty}\ge \|c\|_{1,\infty}$.
  
  By Lemma \ref{lem:weakbound}, since we know that for all $k$,
  $\limsup_{n\to\infty}h(n)/n^k=\infty$, it follows that
  $\|c\|_{1,\infty}=\infty$.
\end{proof}

The following corollary is closely related to a theorem of Jones and
Wierdl \cite{JonesWierdl} (the hypothesis and conclusion are both
weakened).

\begin{cor}
  If $(a_n)$ is an increasing sequence of real numbers with the
  property that for all $\epsilon>0$, $a_n\le n^\epsilon$ for all
  sufficiently large $n$, then $B_Ng$ fails to satisfy a maximal
  inequality. 
\end{cor}

\begin{proof}
  If $a_n\le n^\epsilon$ for all $n\ge n_0$, then $h(n)\ge
  n^{1/\epsilon}$ for $n\ge n_0^\epsilon$.
\end{proof}

If $S$ is an infinite subset of $\mathbb{N}$, we let $S_N$ denote
$\{n\in S\colon n\le N\}$. For a function $f\in L^1([0,1))$, we
consider the averages $A_Nf(x)=1/|S_N|\sum_{n\in S_N}f(nx)$.

\begin{cor}\label{cor:infinite}
  Let $S$ be an infinite subset of $\mathbb N$.  If $S$ has the
  property that
  \begin{equation*}
    \limsup_{N\to\infty} |S_N|/(\log N)^k=\infty\text{ for all $k$,}
  \end{equation*}
  then there exists $f\in L^1$ such that $\limsup A_Nf(x)=\infty$
  almost everywhere.
\end{cor}

\begin{proof}
  By Fact \ref{fact:div}, the conclusion is equivalent to
  establishing the fact that there is no maximal inequality for the
  averages $A_N$. By Lemma \ref{lem:equiv}, this is equivalent to
  establishing that there is no maximal inequality for the averages
  $B_Ng(y)=(1/|S_N|)\sum_{t\in \log(S_N)}g(y-t)$. Since 
  the number of elements of $\log(S)$ up to $K$ is equal to the
  number of elements of $S$ up to $e^K$, which by hypothesis is not
  bounded by any power of $K$, Theorem \ref{thm:growth} gives the
  desired conclusion.
\end{proof}

\begin{thm}\label{th:dichotomy}
  If $S$ is a multiplicative subsemigroup of the positive integers,
  then there is pointwise convergence of $A_Nf(x)$ to $\int f$ for all
  $f\in L^1$ if and only if $S$ is contained in a finitely generated
  semigroup.
\end{thm}

\begin{proof}
  
  Suppose that $\log b_1,\ldots, \log b_d$ are rationally independent.
  Let $B$ be the largest of the $b$'s. It follows that for any $n\ge
  1$, there are at least $Cn^d$ terms of $S\cap [1,B^n]$.
  
  If $S$ is not contained in any finitely generated semigroup, it
  follows that for any $k$, there exist elements $b_1, b_2,\ldots,b_k$
  of $S$ whose logarithms are rationally independent so that the
  hypothesis of Corollary \ref{cor:infinite} is satisfied showing that
  there exists an $f\in L^1$ such that $\limsup A_Nf(x)$ is infinite
  almost everywhere.

  In the case where $S$ is contained in a finitely generated
  semigroup, we make use of an ergodic theorem for amenable group
  actions due to Ornstein and Weiss \cite{OrnsteinWeiss}. It is
  sufficient to establish that the sets $S_N$ defined above form a F\o
  lner sequence. For convenience, we use additive notation.
  Specifically, since by assumption, $S$ is contained in a finitely
  generated semigroup of the positive integers, let the primes that
  appear as factors of elements of $S$ be $p_1,\ldots, p_k$. Given
  $n\in S$, write $n=p_1^{\alpha_1}\cdots p_k^{\alpha_k}$ and we will
  associate $n$ with the vector $(\alpha_1,\ldots,\alpha_k)\in \mathbb
  Z_+^k$. These vectors span a lattice in $\Z^k$ whose dimension we
  will call $d$. Let $L_+$ be the intersection of $\mathbb Z^k_+$ with
  the lattice spanned by the vectors in $S$. In this notation, $S_N$
  corresponds to $\{(\alpha_1,\ldots,\alpha_k)\in S\colon
  \sum\alpha_i\log p_i\le \log N\}$. Clearly the $S_N$ are nested.  It
  remains to establish the following two conditions.

  \begin{align}
    &\text{ For all $n\in
      S$, }\lim_{N\to\infty} |(n+S_N)\bigtriangleup S_N|/|S_N|=0
    \label{Foln1}\\
    &\text{There exists an $M$ such that for all $N$, }|S_N-S_N|\le
    M|S_N|.\label{Foln2}
  \end{align}
  
  The second of these is seen as follows: If $x\in S_N-S_N$, then $x$
  may be expressed as
  $(x_1,\ldots,x_k)=(\alpha_1,\ldots,\alpha_k)-(\beta_1,\ldots,\beta_k)$,
  where $(\alpha_1,\ldots,\alpha_k)$ and $(\beta_1,\ldots,\beta_k)$
  are in $S_N$. It follows that $\alpha_i\le (\log N)/(\log p_i)$ so
  that $|x_i|\le (\log N)/(\log p_i)$. Clearly the number of such
  elements $x$ is bounded above by an expression of the form $C(\log
  N)^d$. On the other hand, by the argument at the start of the proof,
  there are at least $C'(\log N)^d$ elements in $S_N$ so condition
  \eqref{Foln2} holds.
  
  To establish condition \eqref{Foln1}, let $x\in S$. We need to
  estimate the cardinality of $(S_N+x)\bigtriangleup S_N$. Clearly
  this is twice the cardinality of $(S_N+x)\setminus S_N$. This
  difference is contained in $\{(\alpha_1,\ldots,\alpha_k)\in L_+
  \colon \log N<\sum \alpha_i\log p_i\le \log N+\sum x_i\log p_i\}$.
  To estimate this, we will use a crude estimate for the number $L(y)$
  of lattice points in $\{(x_1,\ldots,x_k)\in L_+\colon
  \sum\alpha_i\log p_i\le y\}$. Let $V$ be the $d$-dimensional vector
  space spanned by $S$ equipped with the inherited $d$-dimensional
  Lebesgue measure $\lambda$. Let $F$ denote a convex fundamental
  domain for the lattice $L$ inside the vector space $V$ and let $T$
  denote the set $V\cap \{(\alpha_1,\ldots,\alpha_k)\colon \sum
  \alpha_i\log p_i\le 1\}$. We claim that
  $L(y)=y^d\lambda(T)/\lambda(F)+o(y^d)$. To see this, note that if
  $h$ is the diameter of $F$ and letting $P_\text{int}(y)$ denote the
  set of lattice points in $L_+$ whose $h$-neighborhood lies within
  $yT$ and let $P_\text{ext}(y)$ denote the set of lattice points
  whose $h$-neighborhood intersects $yT$. We now have
  $\text{Int}_{2h}(yT)\subset P_\text{int}(y)+F\subset yT \subset
  P_\text{ext}(y)+F \subset B_{2h}(yT)$ so that
  $|P_\text{int}(y)|\lambda(F)\le y^d\lambda(T) \le
  |P_\text{ext}(y)|\lambda(F)$. Clearly we also have
  $|P_\text{int}(y)|\le L(y) \le |P_\text{ext}(y)|$ so it follows that
  $|L(y)-y^d\lambda(T)/\lambda(F)|\le
  P_\text{ext}(y)-P_\text{int}(y)\le \lambda\big(B_h(yT)\setminus
  \text{Int}_h(yT)\big)$. Since this region is contained in the union
  of $d+1$ slabs each of which having bounded thickness and dimensions
  linearly dependent on $y$, this quantity is $O(y^{d-1})$. It now
  follows that $L(\log N+\sum x_i\log p_i)-L(\log N)=O((\log
  N)^{d-1})$ so that $\lim_{N\to\infty}|(S_N+x)\setminus S_N|/|S_N|=0$
  as required.

\end{proof}

\section{Ultimate Badness of Exponential sequences}

\begin{thm}
  For any $k\in\{2,3,\ldots\}$, the sequence $(k^n)$ is ultimately bad
in $L^1$.
\end{thm}

\begin{proof}
  We deny a maximal inequality by carefully using the standard
  lacunarity trick for a rotation of the circle.  Let $(w_t)$ be
  a sequence such that $C_1(w)=\infty$. Let $M$ be a large
  integer and fix a $y\in\N$ such that 
  \begin{equation}
    \label{eq:ychoice}
    m_y=\sum_t [2^{-y}<w_t<2^{t-y}]w_t>M.
  \end{equation}
  Let $n_0=\lfloor 2y\log_k 2\rfloor$ so that $k^{n_0}\approx 2^{2y}$.
  Throughout the proof, $K$ will be used to denote various quantities
  that can be bounded above or below independently of $y$ and $M$.
  (The bounds may however depend on $k$). Let $f$ be the function on the
  circle taking the value $2^y$ on an interval of length
  $3/(k^{n_0}-1)$ starting at $-1/(k^{n_0}-1)$ and extending to
  $2/(k^{n_0}-1)$ and 0 elsewhere so that $\|f\|_1\le K2^{-y}$.

  We now construct a number $\alpha$ such that letting $T$ be the
  rotation of the unit circle by $-\alpha$ and computing the averages
  \begin{equation*}
    A_tf(x)=w_t\frac{1}{2^t}\sum_{n\le 2^t}f(x-k^n\alpha),
  \end{equation*}
  the maximal function $f^*(x)=\sup A_tf(x)$ has weak $L^1$ norm
  greater than $Km_y\|f\|_1$.
  
  Initially divide the circle into intervals of length
  $1/(k^{n_0}-1)$.  These intervals have endpoints whose base $k$
  expansions are periodic with period dividing $n_0$. We label each
  interval by a string of $n_0$ symbols in $\{0,\ldots,k-1\}$ that
  form the repeated block of the left endpoint so that if $B\in
  \{0,1,\ldots,k-1\}^{n_0}$ then $I_B$ is the interval with left
  endpoint equal to $0.\overline B$ in the base $k$ expansion. We
  shall consider only those intervals whose left endpoint's 
  expansion has period exactly $n_0$. This excludes a negligible
  fraction of the intervals.
  
  Consider a $t$ satisfying $8n_02^{-y}<w_t<2^{t-y}$. We say that an
  interval $I_B$ is \textsl{infected} at time $2^t$ if
  \begin{equation*}
    w_t\frac{1}{2^t}\sum_{2^{t-1}<n\le 2^t}f(x-k^n\alpha)>1\text{ for
    all $x\in I_B$.}
  \end{equation*}
  We will show how to choose the digits of $\alpha$'s base $k$
  expansion from the $2^{t-1}$ position to the $2^t$ position in order
  to bound below the number of new intervals infected at time $2^t$.
  Summing these contributions over $t$, will give a lower bound for
  the maximal function as required.
  
  We say that two words are (cyclically) equivalent if one is a cyclic
  permutation of the other. List representatives of all of the
  equivalence classes in some order as $B_1$, $B_2$, $\ldots$.
  Suppose that by time $2^{t-1}$, the intervals corresponding to
  $B_1,\ldots,B_j$ and their cyclic permutations are infected.  We
  then define the binary expansion of $\alpha$ starting from the
  $2^{t-1}$st digit to be concatenations of $B_{j+1}$ until the
  intervals corresponding to the members of the equivalence class
  become infected. At this point, define digits of $\alpha$ to be
  concatenations of $B_{j+2}$ etc. If all of the the equivalence
  classes are exhausted before the $2^t$th digit of the binary
  expansion is defined, this will ensure that the constant in the
  maximal inequality exceeds $K2^y$ which will be sufficient as the
  $y$ can be chosen to be arbitrarily large. We estimate the number of
  intervals that can be infected up to time $2^t$ as follows:
  
  Let $v_r$ be the sequence obtained by cyclically permuting $B_{j+1}$
  to the left $r$ times and let $J_r$ be the interval corresponding to
  $v_r$. We observe that if $n\geq 2^{t-1}+n_0$, we can write $n$ as
  $2^{t-1} +jn_0+r$ for some $j\geq 1$ and $0\le r<n_0$.  In this
  case, for $x\in J_r$, we notice that $f(x-k^n\alpha)=2^y$.  In order
  to be infected, the sum needs to exceed $2^t/w_t$ so we see that
  this needs to be repeated $\lceil 2^{t-y}/w_t\rceil$ times.  After
  this number of repetitions, $\alpha$ starts following the next $B$
  in a similar manner. The number of repetitions of $B$ in each block
  is therefore bounded above by $2^{t-y+2}/w_t$ (the extra 1 being an
  overestimate coming from the fact that we have no control of the
  location of $x-k^n\alpha$ while $j=0$). Since each repetition has
  length $n_0$, the length of the block is bounded above by
  $2^{t-y+2}n_0/w_t$. Since we have $2^{t-1}$ digits available to
  define, we are able to infect the intervals in at least $K\lfloor
  2^{t-1}/(2^{t-y+2}n_0/w_t)\rfloor$ equivalence classes. Since we
  ensured that $w_t>8n_02^{-y}$, we see that the quantity being
  rounded is greater than 1. As each equivalence class has $n_0$
  members, we see that the number of intervals infected is given by
  $K2^yw_t$. The measure of the infected intervals then exceeds
  $K2^{-y}w_t\ge Kw_t\|f\|_1$.  The constant in the maximal
  inequality therefore exceeds $Kl_y$ where
  \begin{equation*}
    l_y=\sum_t[8n_02^{-y}<w_t<2^{t-y}]w_t
  \end{equation*}
  We complete the proof by demonstrating that $l_y+l_{2y}>m_y-9$. 

  We have
  \begin{equation*}
    l_y+l_{2y}\ge \sum_t[2^{-y}<w_t<2^{t-2y}\text{ or }
      8n_02^{-y}<w_t<2^{t-y}]w_t.
  \end{equation*}
  This yields 
  \begin{align*}
    m_y-(l_y+l_{2y})&\le \sum_t[2^{t-2y}\le w_t\le 8n_02^{-y}]w_t\\
    &\le 16y2^{-y}\#\{t\colon 2^{t-2y}\le 8n_02^{-y}\}\\
    &\le16y2^{-y}\#\{t\colon 2^t\le 2^{2y}\}\le 32y^22^{-y}<9
  \end{align*}
  as required.
\end{proof}

\begin{thm}
For $k\in \{2,3,\ldots\}$, the sequence $(k^n)$ is ultimately bad for
$L^p$ when $p>1$.
\end{thm}

\begin{proof}
  We will use condition (3) established in Theorem
  \ref{thm:BPrinc} for ultimate badness. We deal with the case $k>2$. The
  fact that $(2^n)$ is ultimately bad follows from the fact that
  $(4^n)$ is ultimately bad using Remark \ref{rem:subseq}.
  
  For a given subset $J$ of the positive integers, we construct a
  characteristic function $f=1_B$ on $\mathbb Z$ such that $A_jf$ (the
  average over the $j$th dyadic block) takes a value of order 1 on a
  set of size approximately $|B|$, but that $A_jf$ and $A_{j'}f$ are
  disjointly supported for distinct $j,j'\in J$.
  
  Let $J$ be a finite set of integers. We will assume that $J$
  contains no two consecutive integers.  For $j\in J$, let $B_j$ denote
  $(k-1)\cdot k^{2^j}\cdot \{1,2,3,4,\ldots,k^{2^{j+1}}\}$. Let $C_j$
  denote the truncated version $(k-1)\cdot
  k^{2^j}\cdot\{1,2,3,4,\ldots,k^{2^{j+1}}-k^{2^{j}}\}$.
  
  Let $B=\sum _{j\in J} B_j$, $C=\sum_{j\in J}C_j$; and let $f\in l^p$
  be the characteristic function of $B$. By the requirement
  that $J$ contains no two consecutive integers, it follows that each
  element of $B$ may be expressed in only one way as the sum of
  elements of the $B_j$'s. Note that $|C_j|/|B_j|=1-k^{-2^j}$ so that
  since $|C|=\prod_{j\in J}|C_j|$ and $|B|=\prod_{j\in J}|B_j|$, we
  have $|C|\ge |B|/2$.
  
  Let $x\in \mathbb Z$ satisfy $x+k^{2^{j_0}}\in C$, for some $j_0\in
  J$. Let $n=2^{j_0}$.  Let $m\in [2^{j_0},2^{j_0+1})$.  We have
  $x+k^n=\sum _{j\in J} c_j$, where $c_j\in C_j$. In other words,
  $x+k^n=\sum _{j\in J\setminus \{j_0\}}c_j + (k-1)\cdot k^n a$, where
  $0<a\le k^{2n}-k^n$. We now have $x+k^m = x+k^n+(k^m-k^n)$.  Then
  $k^m-k^n=k^n(k^{m-n}-1)= (k-1)\cdot k^n((k^{m-n}-1)/(k-1))$ so that
  $x+k^m=\sum _{j\in J- \{j_0\}}c_j + (k-1)\cdot k^n
  (a+(k^{m-n}-1)/(k-1))$.  Since $a\le k^{2n}-k^n$ and $m-n<n$, it
  follows that $a+(k^{m-n}-1)/(k-1)\le k^{2n}$ ensuring that $x+k^m\in
  B$.  This establishes that for $x\in C-k^{2^{j_0}}$, $A_{j_0}f(x)\ge
  1/2$.
  
  We now show that these sets are disjoint. Suppose that $x$ lies in
  $C-k^{2^l}$ and $C-k^{2^m}$ for distinct $l>m\in J$.  Then we see
  that $k^{2^l}-k^{2^m}\in B-B$.  We show that this gives rise to a
  contradiction as follows.  Note that $B-B= \sum _{j\in J} S_j$,
  where $S_j=(k-1)\cdot k^{2^j}\cdot \{t\colon |t|<k^{2^{j+1}}\}$.
  If $a\in B-B$ has its largest non-zero summand in the $S_j$ block,
  then we see that $k^{2^{j+2}}/2>(k-1)(k^{3\cdot 2^j}+k^{3\cdot
  2^{j-2}}+\ldots)\ge |a|\ge (k-1)(k^{2^j}-k^{3\cdot 2^{j-2}}-k^{3\cdot 
  2^{j-4}}\ldots)>k^{2^j}$ If $j\ge l$, we see that $a$ exceeds
  $k^{2^l}-k^{2^m}$. If $j<l$ then $j\le l-2$ and we see that $a$ is
  smaller than $k^{2^l}-k^{2^m}$. Note that this is where we made use
  of the assumption that $k>2$.
  
  It follows that $\|\sup_{j\in J}A_jf\|_{p,\infty}\ge
  (|J||C|)^{1/p}/2 \ge |J|^{1/p}\|f\|_p/4$.  A standard argument using
  Rokhlin towers similar to (but simpler than) Lemma \ref{lem:equiv}
  allows this example to be transferred to an arbitrary aperiodic
  system.
\end{proof}

\section{Questions and Remarks}
\begin{rem}\label{rem:LpubnotL1ub}
  The sequence of times $a_n=\lfloor \log n\rfloor$ is ultimately bad
  in $L^p$ for all $1<p<\infty$, but not in $L^1$. To see that the
  sequence is ultimately bad in $L^p$, using Theorem \ref{th:lotech},
  we verify that for $f(x)=2\cdot \mathbf 1_{[0,2\log 2]}$ and
  $J\subset \mathbb N$, for $x\in [j,j+1]$, we see that $A_jf(x)\ge
  1$, verifying condition \eqref{it:lotweak} of Theorem
  \ref{th:lotech}.

  To see that the sequence is not ultimately bad in $L^1$, let
  $w_t=1/t$ and note that $C_1(w_t)=\infty$ but $\sup_tw_tA_tf$ is
  bounded above by the regular ergodic maximal function of $f$, which
  has weak $L^1$ norm bounded above by $\|f\|_1$.
\end{rem}

\begin{quest}\label{quest:KhintchL1}
  Is the sequence of operators $M_nf(x)=f(nx\bmod 1)$ ultimately bad
  in $L^1$?
\end{quest}

\begin{rem}
  We remark that Theorem \ref{th:Khintchine} also shows that if
  $s(n)/\log n\to 0$, then there exists an $f\in L^1$ for which
  $1/(ns(n))\sum_{k\le n}f(kx)$ diverges almost everywhere. We also
  pose the following weakening of Question \ref{quest:KhintchL1}.
\end{rem}

\begin{quest}
  Does there exist $f\in L^1([0,1))$, such that $1/(n\log n)\sum_{k\le
  n}f(kx)$ diverges almost everywhere?
\end{quest}

\begin{quest}
  Do there exist lacunary sequences that are not ultimately bad in
  some $L^p$?
\end{quest}

\begin{quest}
  If the sequence of times $(a_n)$ is ultimately bad for $L^1$, does
  it follow that it is ultimately bad for $L^p$ $(p>1)$.
  Remark \ref{rem:LpubnotL1ub} shows that the converse is false.
\end{quest}

\end{document}